# A binomial formula for evaluating integrals

## Khristo N. Boyadzhiev

Department of mathematics

Ohio Northern University

Ada, Ohio 45810, USA

[k-boyadzhiev@onu.edu](mailto:k-boyadzhiev@onu.edu)

## Abstract

In this paper we present a special formula for transforming integrals to series. The resulting series involves binomial transforms with the Taylor coefficients of the integrand. Five applications are provided for evaluating challenging integrals.

Keywords: Binomial transform, binomial identities, integral evaluation, harmonic numbers, dilogarithm, trilogarithm.

Mathematics Subject Classification: 05A19, 26A42, 40A30.

## 1. Main Theorem

We will present a rule for evaluating integrals in terms of series with binomial expressions.

**Theorem 1.1.** *Let* $f(x)$ *be a function defined and integrable on* $(-r, \lambda]$ *for some* $r > 0, \lambda > 0$.

*Let also* $f(x)$ *be analytic in a neighborhood of the origin with Taylor series* $f(x) = \sum_{n=0}^{\infty} a_n x^n$.

*Then we have*

$$\int_0^{\lambda} f(x)dx = \sum_{n=0}^{\infty} \left(\frac{\lambda}{\lambda+1}\right)^{n+1} \frac{1}{(n+1)} \sum_{m=0}^{n} b_m = \sum_{n=0}^{\infty} \left(\frac{\lambda}{\lambda+1}\right)^{n+1} \sum_{k=0}^{n} \binom{n}{k} \frac{a_k}{k+1}$$

*where the sequence* $b_n$ *is the binomial transform of the sequence* $a_n$

$$b_n = \sum_{k=0}^{n} \binom{n}{k} a_k.$$

*In particular, for* $\lambda = 1$ *we have*

$$\int_0^1 f(x)dx = \sum_{n=0}^{\infty} \frac{1}{2^{n+1}(n+1)} \sum_{m=0}^{n} b_m = \sum_{n=0}^{\infty} \frac{1}{2^{n+1}} \sum_{k=0}^{n} \binom{n}{k} \frac{a_k}{k+1}$$

*and for* $\lambda \to \infty$

$$\int_0^\infty f(x)dx = \sum_{n=0}^\infty \frac{1}{(n+1)} \sum_{m=0}^n b_m = \sum_{n=0}^\infty \sum_{k=0}^n \binom{n}{k} \frac{a_k}{k+1}.$$

**Proof.** With the substitution $x = \dfrac{t}{1-t}, t = \dfrac{x}{x+1}$ we get

$$\int_0^\lambda f(x)dx = \int_0^{\lambda/(\lambda+1)} \frac{1}{(1-t)^2} f\left(\frac{t}{1-t}\right) dt = \int_0^{\lambda/(\lambda+1)} \frac{1}{1-t} \left\{ \frac{1}{1-t} f\left(\frac{t}{1-t}\right) \right\} dt$$

$$= \int_0^{\lambda/(\lambda+1)} \frac{1}{1-t} \left\{ \sum_{n=0}^\infty t^n \sum_{k=0}^n \binom{n}{k} a_k \right\} dt = \int_0^{\lambda/(\lambda+1)} \frac{1}{1-t} \left\{ \sum_{n=0}^\infty b_n t^n \right\} dt$$

by using Euler's series transformation formula

$$\frac{1}{1-t} f\left(\frac{t}{1-t}\right) = \sum_{n=0}^\infty t^n \left\{ \sum_{k=0}^n \binom{n}{k} a_k \right\} = \sum_{n=0}^\infty b_n t^n$$

where the sequence $\{b_n\}$ is the binomial transform of the sequence $\{a_n\}$ as described above..

Expanding $(1-t)^{-1}$ as geometric series and using Cauchy's rule for multiplication of two power series we write

$$\int_0^{\lambda/(\lambda+1)} \frac{1}{1-t} \left\{ \sum_{n=0}^\infty b_n t^n \right\} dt = \int_0^{\lambda/(\lambda+1)} \sum_{n=0}^\infty \left\{ \sum_{k=0}^n b_k \right\} t^n dt = \sum_{n=0}^\infty \left(\frac{\lambda}{\lambda+1}\right)^{n+1} \frac{1}{n+1} \sum_{k=0}^n b_k$$

by the property

$$\sum_{k=0}^n \binom{n}{k} \frac{a_k}{k+1} = \frac{1}{n+1} \sum_{k=0}^n b_k$$

(see [1, p. 61]). The interchange of integration and summation is justifies as we work with power series. This way the first part of the theorem is proved. The proof of the second part repeats the same steps.

Differentiating in the theorem with respect to $\lambda$ we come to the following result.

**Corollary 1.2.** *Under the conditions of the theorem we have the representation*

$$f(\lambda) = \frac{1}{(\lambda+1)^2} \sum_{n=0}^\infty \left(\frac{\lambda}{\lambda+1}\right)^n \sum_{m=0}^n b_m.$$

## 2. Applications

Here we give some applications of our theorem in the form of examples.

**Example 1.** In our first example we will evaluate the integral

$$\int_0^\infty \frac{\log(1+t)}{t(1+t)}\,dt\,.$$

We start from the well-known series

$$\sum_{n=1}^\infty H_n\,t^n = \frac{-\log(1-t)}{1-t}$$

(Here $H_n = 1 + \frac{1}{2} + \ldots + \frac{1}{n}$, $H_0 = 0$ are the harmonic numbers). Replacing $t$ by $-t$ and dividing both sides by $t$ we get

$$\frac{\log(1+t)}{t(1+t)} = \sum_{k=0}^\infty (-1)^k H_{k+1}\,t^k$$

and we take $a_k = (-1)^k H_{k+1}$. Equation (9.32) in [1] says that

$$\sum_{k=0}^n \binom{n}{k}\frac{a_k}{k+1} = \sum_{k=0}^n \binom{n}{k}\frac{(-1)^k H_{k+1}}{k+1} = \frac{1}{(n+1)^2}\,.$$

This way

$$\int_0^\infty \frac{\log(1+t)}{t(1+t)}\,dt = \sum_{n=0}^\infty \sum_{k=0}^n \binom{n}{k}\frac{a_k}{k+1} = \sum_{n=0}^\infty \frac{1}{(n+1)^2} = \frac{\pi^2}{6}$$

(see Entry 4.291.12 in [3]).

**Example 2.** Here we will evaluate the difficult integral

$$\int_0^\infty \left(\frac{\log(1+t)}{t}\right)^2 dt$$

The integral can be reduced to the previous one using integration by parts, but we will do it independently for illustrating the method.

We start from the well-known power series [4, (5.5.28)]

$$\frac{\log^2(1-t)}{2t} = \sum_{n=1}^\infty \frac{H_n\,t^n}{n+1}$$

where we replace $t$ by $-t$ and then divide both sides by $t$ to write

$$\frac{\log^2(1+t)}{2t^2} = \sum_{n=1}^\infty \frac{(-1)^{n-1} H_n\,t^{n-1}}{n+1} = \sum_{k=0}^\infty \frac{(-1)^k H_{k+1}\,t^k}{k+2}\,.$$

So we take

$$a_k = \frac{(-1)^k H_{k+1}}{k+2}, \quad \sum_{k=0}^{n} \binom{n}{k} \frac{a_k}{k+1} = \sum_{k=0}^{n} \binom{n}{k} \frac{(-1)^k H_{k+1}}{(k+1)(k+2)}.$$

Now

$$\sum_{k=0}^{n} \binom{n}{k} \frac{(-1)^k H_{k+1}}{(k+1)(k+2)} = \sum_{k=0}^{n} \binom{n}{k} \frac{(-1)^k H_{k+1}}{k+1} - \sum_{k=0}^{n} \binom{n}{k} \frac{(-1)^k H_{k+1}}{k+2}$$

$$= \frac{1}{(n+1)^2} - \sum_{k=0}^{n} \binom{n}{k} \frac{(-1)^k}{k+2} \left( H_k + \frac{1}{k+1} \right)$$

(using again [1, (9.32)]). Next, applying property (5.5) from [1]

$$-\sum_{k=0}^{n} \binom{n}{k} \frac{(-1)^k}{k+2} \left( H_k + \frac{1}{k+1} \right) = -\sum_{k=0}^{n} \binom{n}{k} \frac{(-1)^k H_k}{k+2} - \sum_{k=0}^{n} \binom{n}{k} \frac{(-1)^k}{(k+1)(k+2)}$$

$$= \frac{n + H_n}{(n+1)(n+2)} - \frac{1}{n+2} = \frac{H_n - 1}{(n+1)(n+2)}.$$

It is easy to see that

$$\sum_{n=0}^{\infty} \frac{H_n - 1}{(n+1)(n+2)} = \sum_{n=0}^{\infty} \frac{H_n}{(n+1)(n+2)} - \sum_{n=0}^{\infty} \frac{1}{(n+1)(n+2)} = 1 - 1 = 0$$

and we compute

$$\int_0^{\infty} \frac{\log^2(1+t)}{2t^2} dt = \sum_{n=0}^{\infty} \sum_{k=0}^{n} \binom{n}{k} \frac{a_k}{k+1} = \sum_{n=0}^{\infty} \frac{1}{(n+1)^2} = \frac{\pi^2}{6}.$$

Finally,

$$\int_0^{\infty} \left( \frac{\log(1+t)}{t} \right)^2 dt = \frac{\pi^2}{3}.$$

**Example 3.** Using some well-known generating function we evaluate here the integral

$$\int_0^1 \mathrm{Li}_2 \left( \frac{t}{1+t} \right) dt.$$

Here $\mathrm{Li}_2(x)$ is the dilogarithm [5]

$$\mathrm{Li}_2(x) = \sum_{n=1}^{\infty} \frac{x^n}{n^2} \quad (\,|x| < 1)\,.$$

We have

$$\mathrm{Li}_2\left(\frac{t}{1+t}\right)=\sum_{n=1}^{\infty}\frac{(-1)^{n-1}H_n t^n}{n}\quad(\,|\,t\,|<1)$$

so that we take

$$a_n=\frac{(-1)^{n-1}H_n}{n},\quad\frac{a_k}{k+1}=\frac{(-1)^{k-1}H_k}{k(k+1)}=\frac{(-1)^{k-1}H_k}{k}-\frac{(-1)^{k-1}H_k}{k+1}$$

and using two binomial transform formulas (9.4a) and (9.32) from [1] we have

$$\sum_{k=0}^{n}\binom{n}{k}\frac{a_k}{k+1}=\sum_{k=0}^{n}\binom{n}{k}\frac{(-1)^{k-1}H_k}{k}-\sum_{k=0}^{n}\binom{n}{k}\frac{(-1)^{k-1}H_k}{k+1}=H_n^{(2)}-\frac{H_n}{n+1}.$$

Here

$$H_n^{(2)}=1+\frac{1}{2^2}+...+\frac{1}{n^2},H_0^{(2)}=0\,.$$

This way

$$\int_0^1\mathrm{Li}_2\left(\frac{t}{1+t}\right)dt=\sum_{n=0}^{\infty}\frac{H_n^{(2)}}{2^{n+1}}-\sum_{n=0}^{\infty}\frac{H_n}{2^{n+1}(n+1)}\,.$$

These two series are easy to evaluate. We have (see p. 292 in [2])

$$\sum_{n=0}^{\infty}H_n^{(2)}x^n=\frac{\mathrm{Li}_2(x)}{1-x},\quad\sum_{n=0}^{\infty}\frac{H_n x^n}{n+1}=\frac{\log^2(1-x)}{2x}\quad(\,|\,x\,|<1)$$

and we compute with $x=\frac{1}{2}$

$$\int_0^1\mathrm{Li}_2\left(\frac{t}{1+t}\right)dt=\mathrm{Li}_2\left(\frac{1}{2}\right)-\frac{\log^2(2)}{2}=\frac{\pi^2}{12}-\log^2(2)\,.$$

At the end of the example we used the well-known formula.

$$\mathrm{Li}_2\left(\frac{1}{2}\right)=\frac{\pi^2}{12}-\frac{\log^2 2}{2}\,.$$

**Example 4**. Let $q$ be a positive integer. In this example we will evaluate the integral

$$\int_0^1\frac{x^q}{(1+x)^{q+1}}dx$$

by using our theorem.

Let $a_k=\binom{k}{q}(-1)^k$. Then

$$f(x) = \sum_{k=0}^{\infty} a_k x^k = \sum_{k=0}^{\infty} \binom{k}{q} (-1)^k x^k = \sum_{k=q}^{\infty} \binom{k}{q} (-x)^k = \frac{(-x)^q}{(1+x)^{q+1}} = \frac{(-1)^q x^q}{(1+x)^{q+1}}.$$

$$\int_0^1 f(x) dx = (-1)^q \int_0^1 \frac{x^q}{(1+x)^{q+1}} dx = \sum_{n=0}^{\infty} \frac{1}{2^{n+1}} \sum_{k=0}^{n} \binom{n}{k} \frac{a_k}{k+1}$$

$$= \sum_{n=0}^{\infty} \frac{1}{2^{n+1}} \sum_{k=0}^{n} \binom{n}{k} \binom{k}{q} \frac{(-1)^k}{k+1} = (-1)^q \sum_{n=q}^{\infty} \frac{1}{2^{n+1}(n+1)}$$

$$vvv$$

because ([1, equation 10.28]])

$$\sum_{k=0}^{n} \binom{n}{k} \binom{k}{q} \frac{(-1)^k}{k+1} = \frac{(-1)^q}{n+1}.$$

This way

$$\int_0^1 \frac{x^q}{(1+x)^{q+1}} dx = \sum_{n=0}^{\infty} \frac{1}{2^{n+1}(n+1)} - \sum_{n=0}^{q-1} \frac{1}{2^{n+1}(n+1)} = \log 2 - \sum_{n=1}^{q} \frac{1}{2^n n}.$$

Compare to Entry 3.194.8 in [3].

**Example 5.** In this example we will evaluate the challenging integral

$$\int_0^1 \frac{\log^2(1+x)}{2x} dx.$$

We have the expansion [4, (5.5.2)]

$$\frac{\log^2(1+t)}{2t} = \sum_{k=0}^{\infty} \frac{(-1)^{k-1} H_k t^k}{k+1} \quad (|t| < 1)$$

and here

$$a_k = \frac{(-1)^{k-1} H_k}{k+1}, \quad \sum_{k=0}^{n} \binom{n}{k} \frac{(-1)^{k-1} H_k}{k+1} = \frac{H_n}{n+1}$$

([1, (9.32)]). Also, according to the property [1, (5.7)]

$$\sum_{k=0}^{n} \binom{n}{k} \frac{a_k}{k+1} = \frac{1}{n+1} \sum_{k=0}^{n} \frac{H_k}{k+1}.$$

Therefore,

$$\int_0^1 \frac{\log^2(1+x)}{2x}\,dx = \sum_{n=0}^{\infty} \frac{1}{2^{n+1}(n+1)} \sum_{k=0}^{n} \frac{H_k}{k+1}.$$

A simple computation shows that

$$\sum_{k=0}^{n} \frac{H_k}{k+1} = \frac{1}{2}\left(H_n^2 - H_n^{(2)}\right) + \frac{H_n}{n+1}$$

(an identity interesting by itself). From this

$$\int_0^1 \frac{\log^2(1+x)}{2x}\,dx = \frac{1}{2}\sum_{n=0}^{\infty} \frac{H_n^2 - H_n^{(2)}}{2^{n+1}(n+1)} + \sum_{n=0}^{\infty} \frac{H_n}{2^{n+1}(n+1)^2}.$$

These two sums can be evaluated easily by using the generating functions

$$\sum_{n=0}^{\infty} \frac{(H_n^2 - H_n^{(2)})t^{n+1}}{n+1} = -\frac{1}{3}\log^3(1-t)$$

$$\sum_{n=0}^{\infty} \frac{H_n t^{n+1}}{(n+1)^2} = \frac{1}{2}\log(t)\log^2(1-t) + \log(1-t)\,\mathrm{Li}_2(1-t) - \mathrm{Li}_3(1-t) + \zeta(3)$$

(see [5, p. 303]). Here

$$\mathrm{Li}_3(x) = \sum_{n=1}^{\infty} \frac{x^n}{n^3} \quad (|x| < 1)$$

is the trilogarithm [5]. Setting here $t = \dfrac{1}{2}$ and using the values

$$\mathrm{Li}_2\left(\frac{1}{2}\right) = \frac{\pi^2}{12} - \frac{\log^2 2}{2}, \, \mathrm{Li}_3\left(\frac{1}{2}\right) = \frac{7}{8}\zeta(3) - \frac{\pi^2}{12}\log(2) + \frac{1}{6}\log^2(2)$$

we come to the evaluation

$$\int_0^1 \frac{\log^2(1+x)}{2x}\,dx = \frac{1}{8}\zeta(3).$$